\documentclass [a4paper, twoside, 11pt]{article}
\usepackage[utf8]{inputenc}
\usepackage[T1]{fontenc}
\usepackage[english]{babel}
\usepackage{enumitem}

\usepackage[draft,english]{fixme}
\usepackage{cite}
\usepackage{yfonts}
\usepackage{mathtools}
\usepackage{amsmath, amssymb}
\usepackage{dsfont}
\usepackage[amsmath,thmmarks]{ntheorem} 
\theorembodyfont{\normalfont\slshape}
\theoremseparator{}
\theoremstyle{break}
\newtheorem{thm}{Theorem}[section]

\newtheorem{lemma}[thm]{Lemma}

\theorembodyfont{\normalfont\itshape}

\theoremstyle{nonumberplain}
\theoremheaderfont{%
\normalfont\scshape}
\theorembodyfont{\normalfont}
\theoremsymbol{\ensuremath{\blacksquare}}
\theoremseparator{:}

\usepackage{keyval}
\makeatletter
\newcommand*\DeclareMathSymbolShorthand[2]{
   \begingroup
   \setkeys{DMSS}{name=#2,#1}%
   \if\DMSS@overwrite 
   \else
      \expandafter\@ifdefinable\csname \DMSS@prefix\DMSS@name\endcsname{%
        \def\DMSS@overwrite{00}
      }%
   \fi%
   \if\DMSS@overwrite 
      \expandafter\@firstofone%
   \else\expandafter\@gobble\fi%
   {\protected\expandafter%
        \xdef\csname \DMSS@prefix\DMSS@name \endcsname{%
        \unexpanded\expandafter{\DMSS@format{#2}}}}%
   \endgroup}
\define@key{DMSS}{format}{\let\DMSS@format#1}
\define@key{DMSS}{format*}{\def\DMSS@format{\expandafter#1\@firstofone}}
\define@key{DMSS}{name}{\def\DMSS@name{#1}}
\define@key{DMSS}{prefix}{\def\DMSS@prefix{#1}}
\define@key{DMSS}{overwrite}[true]{%
   \edef\DMSS@overwrite{\csname if#1\endcsname 00\else 01\fi}}
\setkeys{DMSS}{overwrite=false}
\setkeys{DMSS}{format*=}
\newcommand\MakeDeclareMathSetCommand[3]{%
   \expandafter\MakeDeclareShorthandCommandAux\csname math#2format\endcsname
   {#1}{#2}{#3}}
\def\MakeDeclareShorthandCommandAux#1#2#3#4{%
   \newcommand*#1{#4}%
   \newcommand*#2[2][]{%
      \DeclareMathSymbolShorthand{format=#1,prefix=#3,##1}{##2}%
   }}
\MakeDeclareMathSetCommand{\DeclareMathSet}{numbers}{\mathbb}
\makeatother

\makeatletter
\def\blfootnote{\xdef\@thefnmark{}\@footnotetext}
\makeatother

\MakeDeclareMathSetCommand{\DeclareMathNumbers}{numb}{\mathbb}
\MakeDeclareMathSetCommand{\DeclareMathGroup}{group}{\mathrm}
\DeclareMathNumbers{N} 
\DeclareMathNumbers{Z} 
\DeclareMathNumbers{Q} 
\DeclareMathNumbers{R} 
\DeclareMathNumbers{C} 
\DeclareMathGroup{GL} 
\DeclareMathGroup{SL} 
\DeclareMathGroup{Mat} 
\newcommand{\PSL}[0]{\mathbf{PSL}}

\newcommand{\SetH}[0]{\mathds{H}}

\newcommand{\tr}{\text{tr}}
\newcommand{\T}{\mathcal{T}}
\newcommand{\M}{\mathcal{M}}

\renewcommand{\phi}{\varphi}
\renewcommand{\epsilon}{\varepsilon}

\title{The Curvature of the Hitchin Connection}
\author{Jørgen Ellegaard Andersen and Niccolo Skovgård Poulsen}

\begin{document}
\blfootnote{Supported in part by the center of excellence grant "Center for quantum geometry of Moduli Spaces" (DNRF95) from the Danish National Research Foundation. }

\maketitle
\begin{abstract}
  In this paper we calculate the curvature of the Hitchin connection. We further show that a slight (possibly trivial) modification of the Hitchin connection has curvature equal to an explict given multiple of the Weil-Petersen symplectic form on Teichm\"{u}ller space.
\end{abstract}
\begin{center}
\thanks{Dedicated to Nigel Hitchin at the conference {\em Hitchin70},\\ celebrating his 70'th Birthday.}
\end{center}

\section{Indroduction}
In \cite{Hitchin} Hitchin introduce a projectively flat connection in the bundle of quantizations of the moduli spaces $M$ of flat $SU(n)$-connections over a surface of genus $g>1$ with central holonomy around a marked point on the surface. This connection was also constructed in \cite{ADW} by Axelrod, Della Pietra and Witten from a more physical perspective, where it was also establish how it is related to quantum Chern-Simons theory. See also \cite{Andersen2}, where it was shown how these two constructions agree and can be slightly generalised.  Let us here briefly recall the setup. 

The moduli space $M$ is compact and smooth in the co-prime case, i.e. in case when the central holonomy around the special marked point generates the centre of $SU(n)$. In general it has a smooth part $M'$, which consist of the irreducible connections (if $n=2$, then $g>2$ for this to be the case, since $M={\mathbb P}^3$ in the case $(g,n)=(2,2))$. The smooth part $M'$ has a natural symplectic form called  the Seshadri-Atiyah-Bott-Goldmann symplectic form. The Chern-Simons line  bundle ${\mathcal L}$ over $M$ is a prequantum line bundle for $\omega$ \cite{Fr}. By the Narasimhan-Seshadri Theorem \cite{NS1,NS2}, the moduli space further has a natural K\"ahler structure once a complex structure on $\Sigma$ has been choosen. This gives a family of complex structures on the moduli space $M$ parametrized by the Teichm\"uller space of $\Sigma$, which we denote $\T$.  Consider now the trivial $C^\infty(M,\mathcal{L}^k)$-bundle ${\mathcal H}^{(k)}$ over Teichm\"{u}ller space $\T$. Then a Hitchin connection is a conncetion in ${\mathcal H}^{(k)}$, which preserve the sub-bundle of holomorphic sections $H^{(k)}(M_\sigma,\mathcal{L}^k)$, $\sigma\in \T$. Further, we require it is given by adding a differential operator valued one form to the trivial connection in ${\mathcal H}^{(k)}$ 
$$\nabla_V^H = \nabla_V^t+u(V),$$
for all vector fields $V$ on $\T$.
Hitchin found an explicit formula for $u$, which in \cite{Andersen2}  is proven to be given by the following global differential operator
\begin{align*}
  u(V)=\frac{-1}{2n+4k}\left(\Delta_{G(V)}+2\nabla_{G(V)d{F_\sigma}}+4k V'[F]_\sigma\right).
\end{align*}
Here $F_\sigma$ is a Ricci potential for $M_\sigma$ the moduli space with the K\"ahler structure given by the point $\sigma\in\T$. The notation $V'$ indicate that we project $V$ onto the holomorphic directions on $\T$. Finally the symetric two tensor, $G(V)$ is given by $G(V)=V'[g_{M_\sigma}^{-1}]$ and the operator $\Delta_{G(V)}$ is given by
\begin{multline*}
   \Delta_{G(V)}:C^\infty(M,\mathcal{L}^k)\stackrel{\nabla^{1,0}_\sigma}{\xrightarrow{\hspace*{1.4cm}}}C^\infty(M,T_\sigma^*\otimes\mathcal{L}^k)\stackrel{G\otimes Id}{\xrightarrow{\hspace*{1.4cm}}}C^\infty(M,T_\sigma\otimes\mathcal{L}^k)\\ \stackrel{\nabla^{1,0}_\sigma\otimes Id+Id\otimes\nabla_\sigma^{1,0}}{\xrightarrow{\hspace*{3cm}}}C^\infty(M,T_\sigma^*\otimes T_\sigma\otimes\mathcal{L}^k)\stackrel{\tr}{\xrightarrow{\hspace*{1cm}}}C^\infty(M,\mathcal{L}^k)
\end{multline*}
For this Hitchin connection it was shown in \cite{Andersen1}, that the curvature  is given by 
\begin{thm}[{\cite[Theorem 4.8]{Andersen1}}]\label{thm:hitJN}
  The curvature of the Hitchin connection acts by$$F_\nabla^{2,0}=\frac{k}{(2k+2n)^2}P_k(\partial_{\mathcal{T}}c)\quad F_\nabla^{1,1}=\frac{ik}{2k+2n}(\theta-2i\partial_{\mathcal{T}}\bar\partial_{\mathcal{T}}F)\quad F_\nabla^{0,2}=0,$$ on sections of the bundle $H^{(k)}$.
\end{thm}
Here $\theta$ is as defined below in (\ref{eq:theta def}). The one form $c$ on $\T$ with values in $C^\infty(M)$  is given by $$ c(V)=-\Delta_{G(V)}F-dF G(V)dF-2nV'[F].$$ Finally, $P_k(\partial_{\mathcal{T}}c(V,W))$ is the prequantum operator associated with the function $\partial_{\mathcal{T}}c(V,W)\in C^\infty(M)$
\begin{align*}
  P_k(\partial_{\mathcal{T}}c(V,W))=\frac{i}{k}\nabla_{X_{\partial_{\mathcal{T}}c(V,W)}}+\partial_{\mathcal{T}}c(V,W),
\end{align*}
where $X_{\partial_{\mathcal{T}}c(V,W)}$ is the Hamiltonian vector field of the function $\partial_{\mathcal{T}}c(V,W)$. In fact it was observe in \cite{Andersen1}, that since the curvature must preserve the holomorphic sections $X_{\partial_{\mathcal{T}}c(V,W)}=0$ and so $d_M(\partial_{\mathcal{T}}c(V,W))=0$. 

The form $\theta$ is given as follows
 \begin{equation}\label{eq:theta def}\theta({\mu_1},\bar{\mu_2})=\frac{1}{4}g_{\mathcal{M}_{VB}^{n,k}}(G({\mu_1})\omega_{\mathcal{M}_{VB}^{n,k}}\bar G(\bar{\mu_2})),
\end{equation}

In this paper we show that 
\begin{lemma}
\label{lem:hit11}
\begin{equation*}
  F_{\nabla^H}^{1,1}=\frac{ik}{2k+2n}(\theta-2i\partial_{\mathcal{T}}\bar\partial_{\mathcal{T}}F)=-\frac{ik (n^2-1)}{12(k+n)\pi}\omega_{\mathcal{T}}
\end{equation*}
\end{lemma}

And using this we can find a $1$-form $\tilde c$ on $\T$ and we consider 
$$\tilde\nabla^H=\nabla^H+\tilde c\otimes \text{Id}_{H^{(k)}}.$$ 
We remark that that it might be that $\tilde c$ is zero. In any case after this (possible trivial) modification, we can prove that
\begin{thm}\label{thm:hit20}
  The connection $\tilde\nabla^H$  is still a Hitchin connection and has pure $(1,1)$ curvature given by 
  $$F_{\tilde \nabla^H} = \frac{ik (n^2-1)}{12(k+n)\pi}\omega_{\mathcal{T}}.$$
\end{thm}

In section \ref{VB} we briefly recall our K\"{a}hler coordinate construction on the universal moduli space of vector bundles from \cite{AP2}. In the following section \ref{sec:hit11}, we compute the $(1,1)$ part of the curvature of the Hitchin connection using the results of \cite{AP2}. In final section \ref{2,0} we modify the Hitchin connection by adding to it a scalar valued one-form  on Teichm\"{u}ller space tensor the identity of $H^{(k)}$, such that the resulting connection has only curvature of type $(1,1)$.

\section{The Moduli Space of Vector Bundles}\label{VB}

In order to compute the curvature of the Hitchin connection, we will use the local coordinates of \cite{AP2}, which we will now briefly recall. Let $\Sigma$ be a surface of genus two or greater.  Pick a point in $\T\times M$, that is a Riemann surface $X$ and a holomorphic vector bundle $E$ over it. For an element 
$$\mu\oplus\nu\in H^1(X,TX)\oplus H^1(X,{\rm End}E)$$ 
define a map $$\chi^{\mu\oplus\nu}:\SetH\times \mathbf{SL}(n,\numbC)\to \SetH\times \mathbf{SL}(n,\numbC)$$ which is annihilated by the following differential operator
\begin{equation*}
   \bar\partial_\SetH \chi^{\mu\oplus\nu}=(\mu-\frac{1}{2}\tilde g^{-1}_X\tr\nu\otimes\nu)\cdot\partial_\SetH\chi^{\mu\oplus\nu}+\partial_{\mathbf{SL}(n,\numbC)}\chi^{\mu\oplus\nu}\cdot\nu.
\end{equation*}
We will denote the projection to $\SetH$ by $\chi_1^{\mu\oplus\nu}$ and the projection to $\mathbf{SL}(n,\numbC)$ by $\chi_2^{\mu\oplus\nu}$.

The near by points contained in the coordinate neighbourhood in $\T\times M$ are represented by a pair of equivalence classes of representations into $\PSL(2,\numbR)$ and $\mathbf{SU}(n)$ respectively. Let's say our base point corresponds to $\rho_\SetH:\pi_1(\Sigma)\to \PSL(2,\numbR)$ and $\rho_E:\pi_1(\Sigma-p)\to \mathbf{SU}(n)$. Then the point corresponding to $\mu\oplus\nu$ is $$ (\rho^{\mu\oplus\nu}_\SetH,\rho_E^{\mu\oplus\nu})(\gamma)=(\chi_1^{\mu\oplus\nu}(\rho_\SetH(\gamma)(\chi_1^{\mu\oplus\nu})^{-1}(z)),\chi_2^{\mu\oplus\nu}(\gamma z,e)\rho_E(\gamma)(\chi^{\mu\oplus\nu}_2(z,e))^{-1}).$$
We proved in \cite{AP2} that this construction gives coordinates and moreover, we provided a Ricci potential for the total space and in particular, we showed in Theorem 4.2 in \cite{AP2}, that for the Ricci potential on $M_\sigma$, which is found in \cite{ZTVB} fulfils 
\begin{lemma}\label{lem:prehit11}
For a pair of vector fields on $\T$ represented by $\mu_1$ and $\bar\mu_2$ we have that
\begin{equation*}
  2\bar\partial_{\T}\partial_{\T} F(\mu_1,\bar\mu_2)= \tr(\mu_1 P^{1,0}_{{\rm End}E}\bar\mu_2P^{0,1}_{{\rm End}E})-i\frac{n^2-1}{6\pi}\omega_{\T}(\mu_1,\bar\mu_2).
\end{equation*}
Where $P^{0,1}_{{\rm End}E}$ (resp. $P^{1,0}_{{\rm End}E}$) is the projection on harmonic $(0,1)$-forms (resp. $(1,0)$-forms) with values in ${\rm End}E$.
\end{lemma}

\section{The $(1,1)$-curvature of the Hitchin Connection}
\label{sec:hit11}

First we calculate $G(V_\mu)$ in coordiantes, here $\mu$ denotes the betrami differential corresponding to $V$ by the Kodaira-Spencer map. We recall from Hitchin \cite{Hitchin} that $G(V_\mu)(\alpha,\beta)=\int_\sigma V_\mu'[-\star_\sigma]\tr\alpha\otimes\beta$. To calculate the variation of $-\star_\sigma$, we need to fix a harmonic 1-form, $\nu$ on $\Sigma$. We split it into $\nu=\nu_1+\bar\nu_2^T$ at a point $X\in \T$ where $\nu_1,\nu_2$ are harmonic $(0,1)$-forms on $X$ with values in ${\rm End}E$. Then we have that at a point $(X_{\mu\oplus 0},E)$, we can use the quasiconformal maps $\chi_1^{\mu\oplus 0}$ to change the complex structure on $X$, so that the complex structure on $X_{\mu\oplus 0}$ is decribed by a quotient construction of $\SetH$ with the standard structure. Then $\nu$ is given by
\begin{align*}
(\chi_1^{\mu\oplus 0})^{-1}_*\nu=&(\nu_1 (\overline{\partial \chi_1^{\mu\oplus 0}})(d\bar z-\mu\frac{{\partial \chi_1^{\mu\oplus 0}}}{\overline{\partial \chi_1^{\mu\oplus 0}}}dz)\\ &+\nu_2 (\partial \chi_1^{\mu\oplus 0})(-\bar\mu\frac{\overline{\partial \chi_1^{\mu\oplus 0}}}{\partial \chi_1^{\mu\oplus 0}}d\bar z+dz) )\circ(\chi_1^{\mu\oplus 0})^{-1}.
\end{align*}
So we can find the harmonic representativ of $\nu$ at $\mu$, which we denote $\nu^{\mu}$, using the projections on harmonic $(1,0)$-forms and $(0,1)$-forms on $X_{\mu\oplus 0}$ with values in End$E$ to obtain that
\begin{align*}
  \nu^{\mu}=&P^{0,1}_{{\rm End}E}((\overline{\partial \chi_1^{\mu\oplus 0}})(\nu_1 d\bar z -\bar\mu\nu_2 d\bar z) )\circ(\chi_1^{\mu\oplus 0})^{-1})
\\ &+P^{1,0}_{{\rm End}E}((\partial \chi_1^{\mu\oplus 0}(\nu_2 dz-\mu\nu_1dz) )\circ(\chi_1^{\mu\oplus 0})^{-1})
\end{align*}
Now $I[\nu]=[-\star\nu^\mu]$ and as is seen in \cite[Lemma 2.15]{Hitchin} we have that $V_\mu(I)[\nu]=[V_\mu(-\star\nu^\mu)]$, since $[-\star V_\mu \nu^\mu]$ is exact. To calculate $V_\mu[-\star\nu^\mu]$, we pull it back to $X_0$ with $\chi_1^{\mu\oplus 0}$ and find that
\begin{align*}
  (\chi_1^{\mu\oplus 0})_*(-\star)\nu^\mu&=iP^{0,1}_{{\rm End}E}((\overline{\partial \chi_1^{\mu\oplus 0}})(\nu_1 d\bar z -\bar\mu\nu_2 d\bar z) )(-\mu ({\partial \chi_1^{\mu\oplus 0}})^{-1} +(\overline{\partial \chi_1^{\mu\oplus 0}})^{-1})
\\ &-iP^{1,0}_{{\rm End}E}((\partial \chi_1^{\mu\oplus 0})(\nu_2 dz+\mu\nu_1dz) )(\bar\mu(\overline{\partial \chi_1^{\mu\oplus 0}})^{-1}+(\partial \chi_1^{\mu\oplus 0})^{-1}).
\end{align*}
When we evaluate this at $\epsilon\mu$ and differentiate with respect to $\epsilon$, then most of the terms have explicit factors of $\epsilon$ and are quickly seen to contribute $-iP^{1,0}_{{\rm End}E}\mu\nu_1-i\mu P^{0,1}_{{\rm End}E}\nu_1$, at $\epsilon =0$. Now the only terms remaining are
 $$P^{0,1}_{{\rm End}E}((\overline{\partial \chi_1^{\mu\oplus 0}})(\nu_1 d\bar z )(\overline{\partial \chi_1^{\mu\oplus 0}})^{-1})$$
 and
 $$P^{1,0}_{{\rm End}E}((\partial \chi_1^{\mu\oplus 0})(\nu_2 dz) )((\partial \chi_1^{\mu\oplus 0})^{-1}).$$ 
 The harmonic projections are given as $P^{0,1}_{{\rm End}E}=I-\bar\partial\Delta_0^{-1}\bar\partial^*$ and $P^{1,0}_{{\rm End}E}=I-\partial\Delta_0^{-1}\partial^*$. When we differentiated these with respect to $\epsilon$ the $I$'s will disappear  and either the first or last $\bar\partial$ or $\bar\partial^*$ (resp. $\partial$ or $\partial^*$) in $\bar\partial\Delta_0^{-1}\bar\partial^*$ (resp. $\partial\Delta_0^{-1}\partial^*$) will not be differentiated. In the first case, we have an exact contribution, which does not change the cohomology class. In the second case the term will be zero, since $\nu\in\ker \bar\partial^*$ ($\bar\nu^T\in\ker\partial^*$). We now conclude that
\begin{equation*}
  V_\mu(I)[\nu]=[-iP^{1,0}_{{\rm End}E}\mu\nu_1-i\mu P^{0,1}_{{\rm End}E}\nu_1].
\end{equation*}
And so we must have that 
$$G(V_\mu)(\nu_1,\nu_2)=-2i\int_\Sigma \mu\tr\nu_1\nu_2,$$ and thus 
$$G(V_{\bar\mu})(\bar\nu_1^T,\bar\nu_2^T)=2i\int_\Sigma \bar\mu\tr\bar\nu_1^T\bar\nu_2^T.$$

Now that we have an expression in our coordinates for $G(V_\mu)$ at the center point, we can calculate (\ref{eq:theta def}) in local coordinates
\begin{multline*}
  G(V_{\mu_1})\omega_{\mathcal{M}_{VB}^{n,k}}\bar G(V_{\bar\mu_2})_{i\bar j}\\=\left(\sum_{j,l}-2i\int_X{\mu_1}\tr \bar\nu_i^T\bar\nu_j^T(-I\int\tr\nu_j\wedge\bar\nu_l^T)2i\int_X\bar{\mu_2}\tr \nu_l\nu_k\right).
\end{multline*}
Also recall that at the center point we have chosen our basis of $\nu_i$'s to be orthonormal and so $P^{0,1}\alpha=-i\sum_i \nu_i\int_\Sigma\tr\alpha\wedge\nu_i$ and so we obtain that
\begin{align*}
 G(V_{\mu_1})\omega_{\mathcal{M}_{VB}^{n,k}}\bar G(V_{\bar\mu_2})_{i\bar j}=4i\left(\int_X{\mu_1}\tr \bar\nu_i^TP^{1,0}(\bar\mu_2\nu_j) \right).
\end{align*}
Contract with the metric and using that $\tr P^{0.1}F=\sum_i\int_\Sigma (F\nu_i)\wedge\bar\nu_i^T $, we get that
\begin{align*}
  \theta(\mu_1,\bar\mu_2)=i\tr(\mu_1P^{0,1}\bar\mu_2 P^{1,0})
\end{align*}
Thus by Lemma~\ref{lem:prehit11} and Theorem~\ref{thm:hitJN}, we have proved Lemma~\ref{lem:hit11}.

\section{Modification of $(2,0)$-part of the Curvature}\label{2,0}

In this section we prove Theorem~\ref{thm:hit20}. First we observe that by the result of the previous section we can use the Bianchi idendity for the curvature to conclude that the $(2,0)$-part of the curvature of the Hitchin connection is $\bar\partial_{\T}$ closed, and hence $d_{\T}$ closed by the following argument. We let $V',W'$ be holomorphic vector fields on $\T$ and $U''$ anti-holomorphic. Then the Bianchi identity gives
\begin{equation*}
  0=U'' (F^{2,0}(V',W')-V'(F^{1,1}(W',U''))+W'(F^{1,1}(U'',V')).
\end{equation*}
But since $F^{1,1}$ is proportional to the symplectic form on $\T$, we get that 
$$-V'(F^{1,1}(W',U''))+W'(F^{1,1}(U'',V'))=\partial_{\T}F^{1,1}(V',W',U'')=0.$$ 
We conclude that $\bar\partial_\T F^{2,0}=0$. Finally we recall form \cite{Andersen1} that $d_M F^{2,0}=0$ as well. Now use that $ F^{2,0}$ is mapping class group invariant, so it pushes down to a closed $(2,0)$-form on the moduli space $\M_g$ of genus $g$ curves.

To proceed further we need to assume that $\Sigma$ has genus three or greater, since this assumption will imply that the following two statement are true.
\begin{itemize}
\item \label{zaal}The moduli space of genus $g\geq 3$ curves, $\M_g$, contains complete curves. This means that there exist a complex surface $S$ and a holomorphic embedding $S\to \M_g$.\ For explicit construction see \cite{Zaal} for genus $3$ and for higher genus references there in.
\item \label{harer} The second thing we need is Harer's result \cite{harer}, that for $g\geq3$ the second cellular homology is
 $$H_2(M_g,\numbC)\cong\numbC.$$ 
\end{itemize}
Harer's result implies that $H^2_{dR}(\M_g,\numbC)\cong\numbC$, since it is dual to $H_2(M_g,\numbC)$. We know that the generator must be $\omega_{\T}$, thus in order to prove that $F^{(2,0)}$ is exact, we need to show that it's class is $0$. We can use the Surface $S$, which is a complex embedding submanifold and we can integrate $F^{(2,0)}$ over it and as it is a $(2,0)$-form the result is $0$, at the same time we know that the integral of $\omega_\T$ is non-zero over $S$ and so the cohomology class of $F^{(2,0)}$ is $0$. This means that there exists a $1$-from $\tilde c$ on $\M_g$  such that $F^{(2,0)}=-d_{\M_g} \tilde c$.

Now we can pull back $\tilde c$ to $\T$ and then define a sligtly modified, but still mapping class group invariant Hitchin connection, as discussed in the introduction. We just need to check that it is still a Hitchin connection. By \cite[Lemma 2.2]{Andersen2} it is enough to prove that
\begin{equation*}
  \frac{i}{2}V[I](\nabla^t_{V})^{1,0}s+\nabla^{0,1}_{M_\sigma}(u(V)+\tilde c(V))s=0
\end{equation*}
But since $\nabla^t+u(V)$ is a Hitchin connection, this reduce to showing $\bar\partial_{M_\sigma}\tilde c(W,V)=0$. But that follows from the defining identity, since $d_\T \tilde c=\partial _\T c$ which is a $(2,0)$ from. To calculate the curvature we se that
\begin{equation*}
  F_{\tilde\nabla}(V,W)=[\nabla_V+\tilde c(V),\nabla_W + \tilde c(W)]=[\nabla_V,\nabla_W]+[\tilde c(V),\nabla_W]+[\nabla_V,\tilde c(W)]+[\tilde c(V),\tilde c(W)].
\end{equation*}
The first term is just the curvature calculated in Theorem~\ref{thm:hitJN}. The two next terms only contribute $-W[\tilde c(V)]+V[\tilde c(W)]=d_\T \tilde c(W,V)$, since $\tilde c$ does not depend on where we are in the moduli space of vector bundles and so commute with the differential operator $u$. The last term is also zero, since multiplication by functions commute, hence we conclude that
\begin{equation*}
  F_{\tilde\nabla}(V,W)=\frac{(n^2-1)k}{6\pi(k+n)}\omega_\T(V,W)+F_\nabla^{(2,0)}(V,W)+d_\T\tilde c(V,W)=\frac{(n^2-1)k}{6\pi(k+n)}\omega_\T(V,W)
\end{equation*}
where the last equality follows by the construction of $\tilde c $, since $F_\nabla^{(2,0)}(V,W)=-d_\T\tilde c(V,W)$. This concludes the proof of Theorem~\ref{thm:hit20}.


\begin{thebibliography}{}

\bibitem
{bers}
 L. Ahlfors \& L. Bers (1960).
\newblock Riemann's mapping theorem for variable metrics.
\newblock {\em Annals of Mathematics}, 72(2):pp. 385--404.



\bibitem
{Andersen1}
J.~E. Andersen \& N.~L. Gammelgaard (2011).
\newblock Hitchin's projectively flat connection, {T}oeplitz operators and the
  asymptotic expansion of {TQFT} curve operators.
\newblock In {\em Grassmannians, moduli spaces and vector bundles}, volume~14
  of {\em Clay Math. Proc.}, pages 1--24. Amer. Math. Soc., Providence, RI.
  
  \bibitem
{Andersen2}
J.~E. Andersen  (2012).
\newblock Hitchin's connection, {T}oeplitz operators, and symmetry invariant
  deformation quantization.
\newblock {\em Quantum Topol.}, 3(3-4):293--325.

\bibitem
{JEANLGML}
J.~E. Andersen, N.~L. Gammelgaard \&  M.~R. Lauridsen (2012).
\newblock Hitchin's connection in metaplectic quantization.
\newblock {\em Quantum Topol.}, 3(3-4):327--357.

\bibitem
{AP}
 J.~E. Andersen \& N. S. Poulsen (2016).
\newblock Coordinates for the Universal Moduli Space of Holomorphic Vector Bundles
\newblock  {\em arXiv:1603.00294}.

\bibitem
{AP2}
 J.~E. Andersen \& N. S. Poulsen (2016).
\newblock An explicit Ricci potential for the Universal Moduli Space Vector Bundles
\newblock  {arXiv:1609.xxxx}.

\bibitem{ADW} S.~Axelrod, S.~Della~Pietra, E.~Witten, "Geometric quantization
of Chern Simons gauge theory.", J.Diff.Geom. {\bf 33} (1991) 787--902.

\bibitem{Fr} D.S. Freed, "Classical Chern-Simons Theory, Part 1",
Adv. Math. {\bf 113} (1995), 237--303.


\bibitem{harer}
{Harer, John}(1983),
  \newblock {The second homology group of the mapping class group of an
              orientable surface},
   \newblock {Invent. Math.},{72}(2):{221--239},

\bibitem
{Hitchin}
N.~J. Hitchin(1990).
\newblock Flat connections and geometric quantization.
\newblock {\em Comm. Math. Phys.}, 131(2):347--380.

\bibitem
{MS}
{ V. B. Mehta \& C. S. Seshadri}(1980).
\newblock {Moduli of vector bundles on curves with parabolic structures},
\newblock {\em Math. Ann.},
  248,(3):{205--239}
  
  \bibitem{NS1} M.S. Narasimhan and  C.S. Seshadri,
"Holomorphic vector bundles on a compact Riemann surface",
Math. Ann. {\bf 155} (1964) 69 -- 80.


\bibitem{NS2} M.S. Narasimhan and C.S. Seshadri, "Stable and unitary vector
bundles on a compact Riemann surface", Ann. Math. {\bf  82} (1965)
540 -- 67.


\bibitem
{NandS}
M.~S. Narasimhan, R.~R. Simha, R. Narasimhan \& C.~S. Seshadri (1963).
\newblock {\em Riemann surfaces}, volume~1 of {\em Mathematical Pamphlets}.
\newblock Tata Institute of Fundamental Research, Bombay.

\bibitem
{ZTVB}
 L.~A. Takhtadzhyan \&  P.~G. Zograf(1989).
\newblock The geometry of moduli spaces of vector bundles over a {R}iemann
  surface.
\newblock {\em Izv. Akad. Nauk SSSR Ser. Mat.}, 53(4):753--770, 911.


\bibitem
{ZTPVB}
  L. A. Takhtajan \& P. Zograf (2008).
\newblock The first {C}hern form on moduli of parabolic bundles.
\newblock {\em Math. Ann.}, 341(1):113--135.

 \bibitem
 {ZTpuncRie}
   L. A. Takhtajan \&  P. G. Zograf(1991).
\newblock A local index theorem for families of {$\overline\partial$}-operators on punctured {R}iemann  surfaces and a new {K}\"ahler metric on their moduli spaces.
\newblock {\em Comm. Math. Phys.}, 137(2):{399--426}.

 \bibitem
 {ZTRie}
   L. A. Takhtajan \&  P. G. Zograf (1987).
\newblock A local index theorem for families of {$\overline\partial$}-operators on {R}iemann  surface.
\newblock {\em Usp.Mat. Nauk}, 42(6):{169-190}(in Russian);Russ. Math. Surv. 42(6) 169-190 (1987).

\bibitem
{Zaal}
Chris Zaal (1995)
\newblock {Explicit complete curves in the moduli space of curves of
              genus three},
 \newblock {Geom. Dedicata}, {56}({2}):{185--196},
 

\bibitem
{Wolpert}
  S. A. Wolpert (1986).
\newblock Chern forms and the {R}iemann tensor for the moduli space of
             curves
\newblock  {\em Invent. Math.}, {85(1)}:119--145.

\end{thebibliography}
\end{document}